\documentclass[12pt]{amsart}
\usepackage[T1]{fontenc}
\usepackage[utf8]{inputenc}
\usepackage[main=english,russian]{babel}
\usepackage{amsmath}
\usepackage{url}
\usepackage{graphicx}
\usepackage{xcolor}

\newcommand{\bq}{\begin{quote}}
\newcommand{\eq}{\end{quote}}
\newcommand{\bi}{\begin{itemize}}
\newcommand{\ei}{\end{itemize}}
\newcommand{\rus}[1]{\foreignlanguage{russian}{#1}}

\title[Mathematics education as a political struggle]{Mathematics education policy as a high stakes political struggle: The case of Soviet Russia of the 1930s}
\author[A.V. Borovik]{Alexandre V. Borovik}
\address{Department of Mathematics, University of Manchester, UK}
\author[S.D. Karakozov]{Serguei D. Karakozov}
\address{Moscow State Pedagogical University, Moscow, Russia}
\author[S.A. Polikarpov]{Serguei A. Polikarpov}
\address{Steklov Mathematical Institute of Russian Academy of Sciences, Moscow, Russia}

\begin{document}
\parskip 5pt

\maketitle

\begin{abstract}
This chapter is an introduction to our ongoing more comprehensive work on a critically important period in the history of Russian mathematics education; it provides a glimpse into the socio-political environment in which the famous Soviet tradition of mathematics education was born. The authors are practitioners of mathematics education in two very different countries, England and Russia. We have a chance to see that too many trends and debates in current education policy resemble battles around mathematics education in the 1920s and 1930s Soviet Russia. This is why this period should be revisited and re-analysed, despite quite a considerable amount of previous research \cite{Karp2010}. Our main conclusion: mathematicians, first of all, were fighting for control over selection, education, and career development, of young mathematicians. In the harshest possible political environment, they were taking potentially lethal risks.
\end{abstract}

\section{Introduction}

In the 1930s leading Russian mathematicians became deeply involved in mathematics education policy. Just a few names: Andrei Kolmogorov, Pavel Alexandrov, Boris Delaunay, Lev Schnirelmann, Alexander Gelfond, Lazar Lyusternik, Alexander Khinchin -- they are remembered, 80 years later, as internationally renowned creators of new fields of mathematical research -- and they (and many of their less famous colleagues) were all engaged in political, by their nature, fights around education. This is usually interpreted as an idealistic desire to maintain higher -- and not always realistic -- standards in mathematics education; however, we argue that much more was at stake.

	Using analysis of contemporary sources, including previously unpublished materials from archives, we show how mathematicians (``the young reformers'', as they are frequently called in Russian literature on the history of mathematics) used their deep involvement in mathematics education and education policy to insist that mathematics as an academic discipline had the right to exist. This could not be taken for granted: for example, at that very time genetics had been systematically suppressed \cite{Lyubishchev2006,Strunnikov1989}. Moreover, mathematics was the first of the sciences to be attacked as an idealistic deviation from the Marxists ``dialectical materialism'' -- it happened in 1929--1931. Mathematics survived, perhaps because it was one of the first such attacks, at the time when the machine of repression was not yet working at its full capacity.

	The starting point of mathematicians' argumentation was an ideological thesis:
\bi
\item The new socialist society (being built in harsh socio-economic conditions) needed a step improvement in mathematics education.
\ei

We show that they used this to prove two points:
\bi
\item academic mathematicians were indispensable as guardians of rigour and the quality of mathematics education;
\item the health of mathematics as a profession required a steady supply of bright young people with good school background in mathematics.
\ei
    	
	We argue that there was an undisclosed aim of these mathematicians' political struggle: namely, an attempt to protect mathematics from the same ideological and political pressure which authorities applied to genetics (and, wider, biology), which eventually led to genetics being forbidden as ``bourgeois pseudo-science'', and to the imprisonment and death of leading Russian geneticists.\footnote{Lyubishchev (2006), Strunnikov (1989), Frolov (1998).}

	In our view, another aim was to prevent the dilution of the academic mathematics community by political appointees; by that time, ``Red Professors'' were already quite prominent in Russian universities and dominated the humanities.

	Finally, in the environment of the Stalinist purges, personal survival was very much at stake.

	Without analysing this page in history, it is difficult to assess important later events and trends in Russian mathematics education. For
example:
\bi
\item The ``new math'' reform in Russian school mathematics education initiated in the 1960-70s by Kolmogorov.
\item The Russian tradition (born in the 1930s and still very valuable) of popularisation of mathematics, outreach activities, mathematical olympiads, and later specialist mathematical schools.
\ei
	In the literature, this theme is covered, in our opinion, inadequately and is occasionally dominated by ideological battles inherited from the past.\footnote{See, for example \cite{Kostenko2013}. }

	A study of this episode in history is also relevant for education policies of today. The authors of this paper are practitioners of mathematics education and education policy and will write on that matter elsewhere.

	A few disclaimers are due.

	This paper is only the start of a big research programme. We expect it to be a hard task.

 	First of all, we are trying to assess something which, for variety of political and cultural reasons, was not properly analysed before. We sit on a mass of material: about 5GB of hundreds of books, papers, archive materials. Only tiny bits of it are reflected in the paper; it simply lists a few key points of our approach. We limit ourselves to highlighting a few pivotal episodes from history, and to introducing people who are likely to play principal roles in any systematic narrative. And we try not to jump to easy conclusions.

	Perhaps the greatest omission from this paper -- again for lack of space -- is discussion of the connection between these struggles in mathematics education with the Luzin affair \cite{Demidov1999,Neretin2020}.

	We hope to publish, in a series of papers, a much more detailed account. This will inevitably involve publication, in some form, of large fragments of Russian language sources and their English translations.

\section{Professor L. L. Mishchenko: Letter to the Academy of Sciences of the USSR}
And now let us look at what may at first sight appear to be a relatively  insignificant document from 1935 -- but the one which is in fact representative of developments. It shows the dramatic human conditions at the heart of our story.

	This document was discovered by chance in archives\footnote{\foreignlanguage{russian}{Архивы Российской Академии Наук, Фонд 389, опись
1, дело 6}, \url{http://isaran.ru/?q=ru/opis&ida=1&guid=C7C839A1-D6F1-2387-A3F2-65FF39F60967},  \foreignlanguage{russian}{Протокол заседания
Комиссии по средней школе при Группы математики АН СССР, переписка с Наркоматом просвещения РСФСР и профессорами Фихтенгольцем Г.М. И
Тартаковским В.А. о преподавании математики в средней школе}.}: a letter with criticism of the standard school geometry textbook, sent in 1935 to the Academy of Sciences of the USSR from a remote settlement in Kazakhstan. The author of the letter gives his name as Professor Leonid Leonidovich Mishchenko. What follows is an English translation of the letter.

\begin{quotation}
To the Division of Physic-Mathematical Sciences of the Academy of Sciences of the USSR

The so-called stabilised\footnote{``Stabilised'' is translation of the Russian \foreignlanguage{russian}{стабильный}; this was an official
term for textbook which were supposed to be used over extended period of time and regularly reprinted.}  textbooks are generally not at the
required level, but textbooks on mathematics stand out as record-holders in exactly the negative sense. Especially illiterate are books by
author collectives 1) Brusilovsky and Gangnus, and 2) Gurwitz and Gangnus.

		I believe that the Academy of Sciences cannot ignore the glaring facts of the official distribution of mathematical ignorance in the USSR. On a practical level, the Academy of Sciences can instruct individuals and commissions to analyse the most dubious textbooks and then make their authoritative decision in a documented report.

		In this letter I offer, as a random illustration to the thought formulated above, one of the problems-theorems given by Gurwitz and Gangnus and analysed in the attached brief note.

		If the Academy of Sciences were to take the trouble to subject the textbooks on mathematics and mechanics to analysis, then, for my part, I would not refuse to take a most active part in this work.

L. Mishchenko

Address: Kazakhstan, post office Shchuch'e, Lunacharsky Str., 25. Professor Mishchenko Leonid Leonidovich.
\end{quotation}

This cover letter was followed with a text with analysis of one of exercise problems in the textbook.

\begin{figure}[h]
\begin{center}
\includegraphics[width=5in]{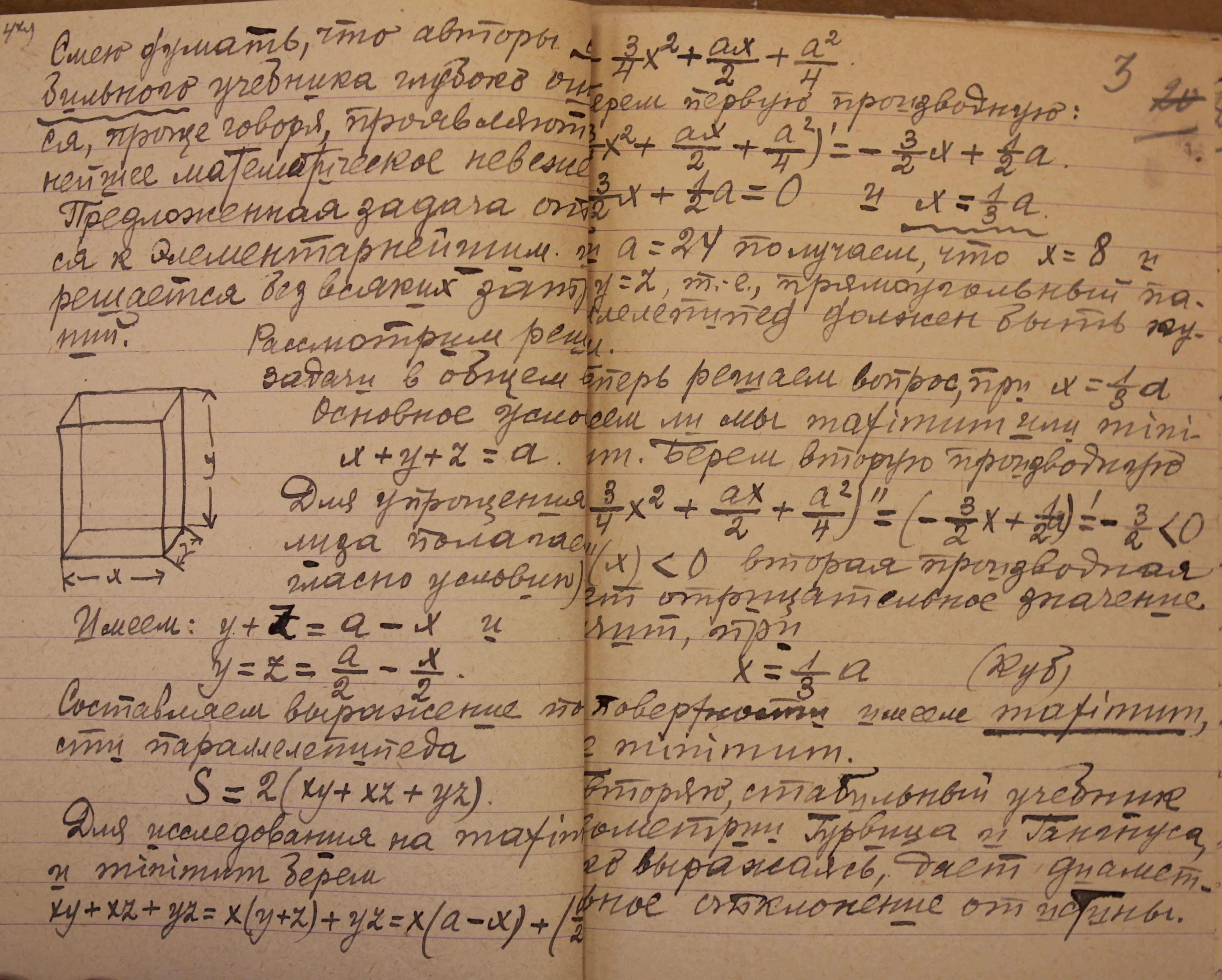}
\end{center}
\caption{A fragment of a letter from L. L. Mishchenko.}
\label{fig:icing}
\end{figure}

Let us take the book: Yu.O. Gurwitz and R.V. Gangnus. A systematic course of geometry. Part 2. Stereometry. Textbook for 8th and 9th forms of secondary school. 1934, Moscow, 2nd Edition, corrected.  Approved by NarKomPros RSFSR.

	On p. 65 in paragraph ``Questions and exercises'' under No. 2 we find the following problem:
Among all rectangular parallelepipeds with the sum of three dimensions equal 24, the cube has the smallest surface area. Check this statement, taking into consideration that [lengths of edges] are integers.
I dare to think that the authors of the stabilised textbook are deeply mistaken, in other words, they show complete mathematical ignorance.

		The proposed problem is one of the most elementary and is solved without any difficulty.
We have to admit that the problem appears to be not very well posed from a pedagogical point of view (even if ``minimum'' is replaced, as it should be, by ``maximum''). Indeed, what were pupils supposed to do? Check every positive integer combination $x + y + z = 24$?

	Mishchenko then gives a standard solution based on differential calculus and the second order derivative test and concludes that at $x = y = z$ (that is, in the cube) the surface has maximum, and not minimum. He continues:
\bq
``I reiterate, the stabilized textbook of geometry by Gurwitz and Gangnus gives a diametrically opposite deviation from the truth.
	I see it as my duty to protest against these outrageous facts of official distribution of mathematical ignorance in RSFSR.
\begin{flushright}
Professor L. L. Mishchenko\\
6 January 1935\\
 Shchuch'e''
\end{flushright}

\eq

\section{Mishchenko's letter: comments}

Mishchenko focuses his assessment of the textbook by Gurwitz and Gangnus on a single problem (most likely, for lack of paper). Mishchenko rightly points out that the cube has not the minimal but the maximal surface area, and, ignoring the assumption that the edges of the parallelepiped have integer lengths, uses calculus to prove that the cube has maximal surface area.  

His calculation leaves an impression that it was done by an engineer, not a mathematician. Indeed, Mishchenko was a military engineer by education and a man of very wide intellectual interests. He taught mathematics and engineering in a variety of institutions of secondary and higher education, was a publisher and editor of journals popularising science and engineering, published books on an incredible range of topics -- from aviation (still fledging, in pre-World War I era) to speech therapy. He was an officer in the Russian army during World War I, was wounded and decorated.

Any normal society at the stage of development when it desperately needed  education would value a competent teacher like Mishchenko. But his letter was handwritten on sheets of paper from a school notebook. At that time (winter 1935), Mishchenko and his wife lived as political exiles in a god-forsaken small settlement on the Karaganda rail line in the brutally cold wind-swept Kazakh steppes. Two years later,  he was arrested and perished in the Gulag. It remains unknown whether his arrest was related to his letter.\footnote{Source with further bibliographic references: \rus{Электронная Библиотека «Научное Наследие России»,} http://library.ruslan.cc/authors/\rus{мищенко-леонид-леонидович/}, accessed 06 January 2020.}

This sets the socio-political background for (not only mathematical) education in the USSR at that time: namely, the acute shortage of educated people -- who, at the same time, were being systematically destroyed. This was the environment in which young reformers who fought for changes in
mathematics education lived and worked.

It is important to emphasise that in his letter, Mishchenko gives a harsh assessment of the textbook by Gurwitz and Gangnus, but makes no political denunciations; his criticism is only directed at mathematical illiteracy.

The sad irony of this story is that Gangnus himself was sent to the Gulag in 1938 (allegedly, for spying on behalf of Latvia), but was lucky -- he was released already in 1943 and exiled to Murom, a provincial town not far from Moscow. He died in 1949 in Moscow.

Meanwhile, Gurwitz continued a successful academic career and was awarded a number of honours. It perhaps helped that he was teaching (and ran a mathematical circle) in a privileged school in Moscow which taught children of the highest tier of the Communist Party elite -- including Stalin's children. He had a reputation for being an excellent teacher as evidenced, for example, by one of his former students \cite{Nikolaev2002}.

\section{Further developments}

It remains unknown whether Mishchenko got any response from the Academy of Sciences. But perhaps his letter contributed to the fate of the textbook by Gurwitz and Gangnus -- it was discovered, unnumbered and not entered in the list of contents, in the same folder in the archives of the Academy of Sciences\footnote{\rus{Архивы Российской Академии Наук, Фонд 389, опись 1, дело 6.}} which contains documents from 1935 and 1936 related to the vicious attack by the young reformers on this particular textbook. One can see how Kolmogorov, Alexandrov, Lyusternik, Delaunay, Schnirelmann, Fichtengolz were appointed to the key committees in the Academy of Sciences and People's Commissariat of Education leading eventually to the following decision of the Mathematics Committee of the People's Commissariat of Education.
\bq
``The textbook of Gurwitz and Gangnus, as a hackwork and an illiterate presentation of the foundations of geometry, has to be immediately withdrawn from schools. The Committee deems the correction of this book in subsequent editions as impossible and suggests, as a temporary measure, publication of one of the geometry textbooks existing in Russian or foreign literature or available in manuscript form. The Committee takes responsibility for recommending and subsequent editing of one of the books.''
\eq
It was followed by publication of papers by Fichtenholtz \cite{Fichtenholtz1937} and  Schnirelman \cite{Schnirelman1937}, evidence of consolidation of the victory of the young reformers and their influence and sometimes control of curriculum development and publication of textbooks for schools. Fichtenholz wrote, among other things, about the need to re-hire experienced ``old'' school teachers.  This would not help Mischenko though: he was arrested on 2 November 1937, on 5 December 1937 sentenced by an infamous Trojka to 10 years in labour camps; his date of death is unknown. It is almost certain that he never learned about the academic debate around the book by Gurwitz and Gangnus.

	The limited space does not allow us to explain the historic background: World War I 1914--18; the Civil War 1918--1921; the ruined economy of Russia; famine; destruction of the old school system replaced by the experimental and vaguely defined concept of ``labour school'' and other bizarre, from the modern point of view, experiments in education in 1918--1930.

	The more pragmatic approach to education and the return to more traditional methods started in 1931\footnote{Marked by decrees of the Central Committee VKP(b) ``On primary and secondary school'' of 25 August 1931 and ``On study programmes and regime in primary and secondary school'' of 25 August 1932. Quite interesting is the directive of Narkompros RSFSR ``On the results of the 1933/34 school year in teaching of principal subjects in primary and secondary school''; it reformulates the Decrees of CC VKP(b) in clear and succinct terms.}, but the new education system (which survived until the collapse of the Soviet Union in 1991) was developing in the environment of the Great Break (the words coined by Stalin in 1929), forced collectivisation and industrialisation, famine in Ukraine, start and escalation of mass repressions, deep distrust of, and systematic repressions against, ``old specialists'' -- this category included essentially all people with pre-revolutionary (that is, before 1917) university education.

	Remarkably, some very high quality popular and semi-popular books on mathematics were published in Leningrad in 1920--30s. But on the national scale, the lack of any system in the production of new school mathematics textbooks is astonishing. For anyone familiar with publication of mathematics books, it is obvious that most textbooks did not pass through a proper refereeing and copy-editing process. Busev \cite{Busev2007,Busev2009} provides a good analysis of reforms of mathematics education and production of textbooks. But the episode with Mishchenko suggests a more serious underlying reason for chaos: the destruction of the mathematics education community. The involvement of professional research mathematicians became really needed -- mathematical education in 1918--30 had lost its error-correcting capacity.

	The sad state of mathematics textbooks produced in the 1930s in the chaos of reforms (and, in particular, of the new ``stable textbooks'') made them an easy target and a welcome chance for the young reformers to prove that mathematicians were indispensable.

	Returning to Gurwitz and Gangnus, Busev \cite{Busev2008} claims with reference to a later highly critical review of the textbook, that about 300 typos and mistakes were identified in the book -- indeed, it was messy. Busev also suggests a plausible explanation for its failure: both Gurwitz and Gangnus were well-educated ``old specialists'', but were conditioned in their involvement in the 1920s in production of study guides for adult education, where mathematical rigour was replaced by explanation of practical examples. He also refers to the ``chaos'' left by the experiments of the 1920s, and to lack of proper copy editing in publishing.

\section{What was at stake?}

We should have no illusion: the very survival of mathematics was at stake. After all, the 1930s were the period of escalating prosecution and step-by-step destruction of genetics -- and geneticists. The peak of the campaign, that is, the VASKhNIL (All-Union Academy of Agricultural Sciences) congress of 1948 was  its most public and triumphal moment; but the systematic suppression started in the 1930s \cite{Strunnikov1989}, and observers from academic circles watched it with increasing alarm. Savina \cite{Savina1995} analyses how these developments were reflected in diaries, letters, and contemporaneous notes by Vernadsky, one of the prominent members of the Academy of Sciences; it is obvious that Lysenko, the principal leader of anti-geneticists, a ``barefoot professor'' as he was described in the official propaganda, was perceived in the academic circles as a very serious danger, the most dangerous breed of a ``Red Professor'', a political appointee to the academia.

	Mathematics was perhaps the first science to be attacked as an idealistic deviation from the Marxists ``dialectical materialism'' -- it happened in 1929--31; see Leifert \cite{Leifert1931} as an example of one of the battles (and political games the mathematicians were forced to play) and Bogolyubov  \cite{Bogolyubov1991} for analysis of the wider landscape. Mathematics survived the ideological attacks of 1929--31, perhaps because at that time the machine of terror was not yet working at its full capacity.

	However, a new large scale attack could start at any time. Freedom -- and life!  -- could not be taken for granted.

	There was also a real threat that the academic mathematical community might be diluted and dominated by ``Red Professors''. By mid-1930s, ``Red Professors'' dominated the humanities -- and were successfully destroying genetics.

 	What would you do if you were a mathematician?

	It is our conjecture (so far supported by all documentary materials that we had a chance to see), that the young reformers started to learn -- and very efficiently -- the rules of the political game.

\section{Meritocratic elitism -- but with political filters}

This subtitle is perhaps the best description of the Soviet system of education as it was shaped in 1930s and existed for several decades.

	The speech by Andrei Bubnov, People's Commissary of Education in 1929--37 at XVII Party Congress of 1934 is quite illuminating. Besides what would now be described as ``widening participation'', ensuring the steady progress of working class children through the school system, he also emphasises a different task which could be formulated as
\bq
Educating the new generation of loyal to the regime high level specialists for the military, industry, science, medicine.
\eq
Some examples given by Bubnov in his speech:
\bq
``Look at Kamai -- professor at Kazan University, Tatar, former docker, now the author of a number of research works in the area of organic compounds of phosphorous and arsenic.'' \cite[p. 114]{VKP(b)1934}
\eq

We should not forget that ``organic compounds of phosphorous and arsenic'' were an obvious euphemism for ``precursors of chemical weapons''. Bubnov's previous post was the Head of the Main Political Directorate of the Red Army. 	The directive of Narkompros RSFSR on the results of the 1933/34 school year contains an interesting line: criticism of teaching geography at school for insufficient attention to giving students practical skills of orientation on the map. It is hard to avoid an impression that Bubnov cared not only about chemical weapons, but also about training of sufficient number artillery officers and air pilots. This was the time of mass conscription armies, and in the hour of need the country had to be prepared to conscript, say, 50 thousand healthy young men with a good knowledge of trigonometry and train them as artillery officers. Indeed, World War II showed that a mathematics graduate could be trained as an artillery officer in a week.

	This explains Bubnov's respect for pure mathematics:
\bq
``Our universities [\dots ] already have young scientists who grew up after the October Revolution -- just look at MGU [Moscow State University] professors Schnirelman\footnote{ Lev Schnirelman: In 1933 he was elected a Corresponding member of the Academy of Sciences of the Soviet Union; committed suicide on 24 September 1938, after being forced to become an informer of NKVD (the secret police). This happened soon after Bubnov was shot as an “enemy of the people” (1 August 1938).} and Gelfond\footnote{Alexander Gelfond: In 1939 he was elected a Corresponding member of the Academy of Sciences of the Soviet Union for his works in the field of cryptography (and not number theory, as one would expect). According to Vladimir Arnold, https://docplayer.ru/40310638-Vladimir-arnold-opasatsya-kompetentnyh-sopernikov-ochen-estestvenno-dlya-nachalnikov.html, during the World War II Gelfond was the Chief Cryptographer of the Soviet Navy.}, already well known for their scientific work in mathematics -- they are young, they were born in 1905.'' \cite[p. 114]{VKP(b)1934}
\eq
It was not number theory (the specialism of Schnirelman and Gelfond) that was needed -- people with mathematics education were needed as a strategic resource of the state. But these mathematicians had to be politically loyal. This opened to the young reformers a window of opportunity.

\section{Mathematical olympiads and ``mathematical outreach''}

The battle for stabilised textbooks was only one of many facets of the work of the young reformers in education. The most interesting and completely new development at that time was the creation of the mathematics olympiad system.

	The first mathematical olympiads were organised in 1934 in Leningrad by Boris Delaunay \cite{Chistyakov1935}, and in 1935 in Moscow by Pavel Alexandrov and other Moscow mathematicians \cite{Bonchkovsky1936}.

	Mathematicians invested remarkable energy and effort in this project. Why?	Because olympiads and other outreach activities were giving them some influence on who is coming into mathematics, and gave a chance to protect university mathematics from political appointees.

	In the Moscow mathematics olympiad of 1935,  Pavel Alexandrov was the chairman of the organising committee. Next year, he wrote in an introduction to a little booklet with problems and solutions of this olympiad:
\bq
``The olympiad is the first entry of future mathematicians into the
mathematical arena. It should help us to select these future
mathematicians from among our youth, it should help us to provide
opportunities for their further mathematical development and education.''\footnote{P.S. Alexandrov's \emph{Foreword} to \cite{Bonchkovsky1936}.}
\eq
Here we see an unashamed elitism in a supposedly egalitarian country -- but, in view of the previous discussion, this is not so surprising. What is really astonishing is the phrase
\bq
``It should help us to select these future mathematicians from among our youth.''
\eq
Selection and promotion of the ``cadre'' was the ultimate monopoly of the ruling Communist Party -- just read a random speech at its XVII Congress \cite{VKP(b)1934}. The young reformers offered, at the right moment, their services to the Party, thus ensuring some degree of their own control over the supply of fresh blood to the top tier of the mathematical profession. In mathematics, ``Red Professors'', recruited from the working class party activists (or Young Communist League activists), could not compete with people who were much more deeply educated and who started their development as mathematicians within the olympiad system.

Mathematicians dared to ask for autonomy in selection and development of their own. It could be seen in other documents of the epoch, for example, in Resolutions of the Second All-Union Congress of Mathematicians which took place in 1934 \cite{VMS1935}: a special resolution was about olympiads (p. 56), and it called for ``identification of gifted youths'', stated that
\bq
``Universities might use olympiads for recruitment of students to mathematical, mechanical and physics departments.''
\eq
and requested funding from the Narkompros (the Ministry of Education) for running olympiads and related activities.

	Apparently, these requests were met by the authorities. Starting from 1935, mathematical olympiads and related activities, first of all, mathematical circles, flourished in Leningrad and Moscow, with \emph{cr\`{e}me de la cr\`{e}me} of research mathematicians actively involved; what is important, new didactical approaches to exposition of non-trivial mathematics were invented and tested. A good and detailed narrative of these remarkable developments can be found in Boltyansky and Yaglom  \cite{Boltyansky1965}.

	This was a two-pronged action of mathematicians aimed at establishing themselves as guardians of standards of mathematics education: one was criticism of school textbooks and involvement in curriculum development, and another one showing by practical work that they can identify and nurture mathematical talent. It was successful and helped mathematicians to improve their positions on the ideological front. We have no space to discuss that in any detail, but can quote a couple of contemporary documents.

	In 1941, Otto Schmidt, a mathematician, a member of the Academy of Sciences, and a high ranking governmental official, was able to formulate the role of mathematicians as specialists who maintain the strategically important mathematical culture of the country:
\bq
``Not only us, professionals of science, but the whole country was happy to learn about solution of a problem set 150 years ago. This problem was solved by academician Vinogradov who proved that every odd number could be written as a sum of three prime numbers. Is this needed at the practical level? No. Maybe it is not needed at all? On the contrary, it is much needed, because the culture, the mathematical culture depends on the level of these works in pure mathematics and theoretical physics. You all know that this is not needed for each of us, but we all are interested in the highest possible level of mathematical culture in our country, because it is important to be able to solve any mathematical problem and for that it is important to be able to solve problems such as Goldbach's problem.

The level of mathematical culture in our country is exceptionally high. One may confidently say that in that respect our countries is on the first and leading place in the world.'' \footnote{	O.Yu. Schmidt. On connection between science and practice, Archive RAS, fund 496, inventory 1, file 248. Quoted from \cite[p. 150]{Dubovitskaya2009}.}
\eq
	Esakov \cite{Esakov1994} gives a tiny, but exceptionally important piece of evidence of Stalin's attitude to mathematics. The text of the speech by Lysenko at the infamous session of the All-Union Academy of Agricultural Sciences in 1948 \cite{VASKhNIL1948} was submitted to Stalin for approval. Stalin underlined Lysenko's statement ``Every science is rooted in class [by its nature]'' and wrote on the margin:
\bq
		``Hah-hah-hah \dots\  And mathematics? And Darwinism?''
\eq
So, by 1948 Stalin (and perhaps even earlier) did not believe in the class nature of mathematics. This was a victory for mathematicians and had a profound effect on the fate of mathematics in the Soviet Union. Not every area of science was so lucky.
	
\section{The reformers' legacy}

Mishchenko's letter also serves as a useful indication of an important part of the reformers' legacy. The  ``olympiad'' or ``outreach'' mathematics which was created by them in Russia at that very time deeply transformed the understanding of elementary mathematics.

	For people influenced by the outreach culture of Russian mathematics, a natural way to approach an extremal geometric problem is to look first at its various degenerate cases, in particular, in the problem of Gurwitz and Gangnus, at the case where the parallelepiped is a long rod with sides $\epsilon$, $\epsilon$, and $24 - 2\epsilon$; its surface area is $2\epsilon^2 + 4\epsilon (24 - 2\epsilon)$ and tends to zero when $\epsilon$ decreases -- the entire calculation could be done in the head. This argument is more appropriate in a discussion of a school textbook. After that, it could be observed that even assuming that all three dimensions are integers, the rod $1 \times 1 \times 22$ has surface area $2 + 4 \times 22 = 90$, which is significantly smaller than the surface area of the cube $8 \times 8 \times 8$ which is $6 \times 8^2 =384$. 

	We think it is important to emphasise that Kolmogorov and his comrades-in-arms did more than critique the state of mathematics education: they had created new educational structures -- such as mathematical circles and mathematics competitions, focused on the mathematics of qualitative analytic thinking. They also created a new, previously never existing, cultural system: the advanced level ``outreach mathematics'', and the community which shared its values.

	By the 1970s, a distinctive ``mathematics outreach'' community had reached a considerable degree of development: a loose informal network of academic mathematicians and school teachers of mathematics involved in organising mathematical competitions, running mathematical circles, summer schools, Sunday schools, distance learning by correspondence schools, undergraduate students who helped to run all these activities -- and, of course, schoolchildren who were enthusiastically taking part in them. The vast majority of professional mathematicians, either in academia or in industry, were to some extent nurtured by this community.

\section{Conclusion}

Borovik \cite[p. 371]{Borovik2016} characterised the current situation in mathematics education as a three-way choice:
\bq
Democratic nations, if they are sufficiently wealthy, have three options:
\bi
\item[(A)] Avoid limiting children's future choices of profession, teach rich mathematics to every child -- and invest serious money into thorough professional education and development of teachers.
\item[(B)] Teach proper mathematics, and from an early age, but only to a selected minority of children. This is a much cheaper option, and it still meets the requirements of industry, defence and security sectors, etc.
\item[(C)] Do not teach proper mathematics at all and depend on other countries for the supply of technology and military protection.
\ei
Which of these options are realistic in a particular country at a given time, and what the choice should be, is for others to decide.
\eq

Russian mathematicians of the 1930s offered, fought a political battle for, and successfully implemented, a solution which appeared to be a version of (B). It is an interesting lesson.

First, their project started in the worst possible economic conditions.
\bi
\item The old system of school mathematics education inherited from the pre-revolutionary times was in ruins.
\item The new system, with a wider social base, had just started to be built from scratch.
\item There was acute shortage of teachers.
\item There were no textbook.
\item In modern European history, only Nazi Germany was able to compete with Stalin's regime in brutality and mass crimes against humanity.
\ei

Its accomplishments became obvious only in the period from the late 1950s to the 1960s, when the Soviet mathematics flourished -- see the previous Section.

	In our next paper (to be published elsewhere) we demonstrate that many elements of this solution naturally fit in a model of mathematics education for the 21st century.

	A lot of interesting historic materials were not even mentioned here; we are preparing a much larger text, to be made public in due course.

\subsection*{Acknowledgments.} The first author thanks Tony Gardiner for helpful advice, and we all thank the reviewers for suggestions which improved the paper. The work on this paper was completed when the first author visited the Nesin Mathematics Village in Turkey\footnote{https://nesinkoyleri.org/en/nesin-villages/}, a unique institution which rightly claims its own place in the history of mathematics education.

	And we wish to express our great respect to Leonid Leonidovich Mishchenko, 1882--19??.

\end{document}